\documentclass[11pt]{article} \newcommand{\nc}{\newcommand}  \nc{\ov}{\over}
\nc{\iy}{\infty} \nc{\x}{\xi} \nc{\inv}{^{-1}} 
\renewcommand{\sp}{\vspace{1ex}} \nc{\be}{\begin{equation}} \nc{\ee}{\end{equation}}
\nc{\ch}{\raisebox{.3ex}{$\chi$}}  \renewcommand{\t}{\tau} \nc{\A}{{\rm Ai}}
\nc{\pl}{\partial} \nc{\de}{\doteq} \nc{\dl}{\delta} \renewcommand{\r}{\rho}
 \nc{\tl}{\tilde} \nc{\sgn}{{\rm sgn}\,}
\nc{\D}{{\cal {D}}\,}
\def\on{\Theta}
\begin{document}

\begin{center}{\bf\large
 A System of Differential Equations for the Airy Process}\end{center}

\sp\begin{center}{{\bf Craig A.~Tracy}\\
{\it Department of Mathematics\\
University of California, Davis, CA 95616\\
email address: tracy@math.ucdavis.edu}}\end{center}
\begin{center}{{\bf Harold Widom}\\
{\it Department of Mathematics\\
University of California, Santa Cruz, CA 95064\\
email address: widom@math.ucsc.edu}}\end{center}

\begin{center}\textbf{Abstract}
\end{center}
The Airy process $\t\rightarrow A_\t$ is characterized by its finite-dimensional
distribution functions 
\[\textrm{Pr}\left(A_{\t_1}< \x_1,\ldots,A_{\t_m}<\x_m\right).\]
For $m=1$  it is
known that $\textrm{Pr}\left(A_\t< \x\right)$ is expressible
in terms of a solution to Painlev\'e II.  We  show that
each finite-dimensional distribution function  is expressible in terms of a solution
to a system of differential equations.
\vskip5pt

\begin{center}{\bf I. Introduction}\end{center}

The Airy process $\t\rightarrow A_\t$, introduced by
Pr\"ahofer and Spohn \cite{PS},  is the limiting stationary process for a certain
$1+1$-dimensional local random growth model called
the polynuclear growth model (PNG).   It is conjectured that
the Airy process is, in fact,  the limiting process for a wide class of 
 random growth models.  (This class is called the $1+1$-dimensional KPZ universality
class in the physics literature~\cite{KPZ}.)  The PNG model is closely related
to the length of the longest increasing subsequence in a random permutation~\cite{AD}.
This fact  together with the result of Baik, Deift and Johansson \cite{BDJ} on
the limiting distribution of the length of the longest increasing subsequence
in a random permutation shows that the distribution function
$\textrm{Pr}\left(A_\t<\x\right)$
equals the limiting distribution function, $F_2(\x)$,
 of the largest eigenvalue in the
Gaussian Unitary Ensemble~\cite{Airy}.  $F_2$   is expressible  
either as a Fredholm determinant of a certain trace-class operator (the Airy kernel)
or in terms of a solution to a nonlinear differential equation (Painlev\'e II).
The finite-dimensional distribution functions
\[ \textrm{Pr}\left(A_{\t_1}<\x_1,\ldots,A_{\t_m}<\x_m\right)\]
are expressible as a Fredholm determinant of a trace-class operator
(the extended Airy kernel)~\cite{J, PS}.
It is natural to conjecture~\cite{J, PS}
 that these distribution functions are also expressible
in terms of a solution to a system of differential equations.  It is this
last conjecture which we prove. For $m=2$ this conjecture was also proved,
in a different form, by Adler and van Moerbeke \cite{adler}.

\begin{center}{\bf II. Statement}\end{center}

The Airy process is characterized by the probabilities
\[{\rm Pr}\;\Big(A_{\t_1}<\x_1,\,\ldots\,,A_{\t_m}<\x_m\Big)=\det\,(I-K),\]
where $K$ is the operator with $m\times m$ matrix kernel having entries
\[K_{ij}(x,y)=L_{ij}(x,y)\,\ch_{(\x_j,\iy)}(y)\]
and
\[L_{ij}(x,y)=\left\{\begin{array}{ll} {\displaystyle{\int_0^\iy}} e^{-z\,(\t_i-\t_j)}\,
\A(x+z)\,\A(y+z)\,dz&{\rm if}\ i\ge j,\\&\\ {\displaystyle -\int_{-\iy}^0}e^{-z\,(\t_i-\t_j)}
\,\A(x+z)\,\A(y+z)\,dz&{\rm if}\ i<j.\end{array}\right.\]
We assume throughout that $\t_1<\cdots<\t_m$, and think of $K$ as acting on the $m$-fold 
direct sum of $L^2(\alpha,\,\iy)$ where $\alpha<\min\,\xi_j$.

To state the result we let $R=K\,(I-K)\inv$ and 
let $A(x)$  denote the \linebreak
$m\times m$ diagonal matrix 
${\rm diag}\,(\A(x))$ and $\ch(x)$ the diagonal matrix
${\rm diag}\,(\ch_j(x))$, where $\ch_j=\ch_{(\x_j,\iy)}$. Then we define the matrix functions
$Q(x)$ and $\tl Q(x)$ by
\[Q=(I-K)\inv A,\ \ \ \tl Q=A\,\ch\,(I-K)\inv\]
(where for $\tl Q$ the operators act on the right).
These and $R(x,y)$ are functions of the $\x_j$ as well as $x$ and $y$. We
define the matrix functions $q,\ \tl q$ and $r$ of the $\x_j$ only by
\[q_{ij}=Q_{ij}(\x_i),\ \ \ \tl q_{ij}=\tl Q_{ij}(\x_j),\ \ \ r_{ij}=R_{ij}(\x_i,\,\x_j).
\footnote{We always interpret $R_{ij}(x,\,\x_j)$ as the limit $R_{ij}(x,\,\x_j+)$. These
quantities are independent of our choice of $\alpha$.}\]
Finally we let $\t$ denote the diagonal matrix ${\rm diag}\,(\t_j)$. 

Our differential operator is $\D=\sum_j\pl_j$, where $\pl_j=\pl/\pl\x_j$, and the system 
of equations is
\begin{eqnarray}
{\cal {D}}^2\,q&=&\x\,q+2\,q\,\on\,\tl q\,q-2\,[\t,\,r]\,q,\label{eq1}\\
{\cal {D}}^2\,\tl q&=&\tl q\,\x+2\,\tl q\,q\,\on\,\tl q-2\,\tl q\,[\t,\,r],\label{eq2}\\
\D r&=&-q\,\on\,\tl q+[\t,\,r].\label{eq3}
\end{eqnarray}
Here the brackets denote commutator, $\x$ denotes the diagonal matrix ${\rm diag}\,(\x_j)$
and $\on$ the matrix with all entries equal to one.

This can be interpreted as a system of ordinary differential equations if 
we replace the variables $\x_1,\,\ldots,\,\x_m$ by
$\x_1+\x,\,\ldots,\,\x_m+\x$, where $\x_1,\,\ldots,\,\x_m$ are fixed and $\x$ variable. 
Then $\D=d/d\x$, and the $\x_j$ are regarded as parameters. 

To get a representation for $\det\,(I-K)$ observe that
\be\pl_j\,K=-L\,\dl_j,\label{plK}\ee
where the last factor denotes multiplication by the diagonal matrix with all entries zero 
except for the $j^{\scriptstyle \rm th}$,
which equals $\dl(x-\x_j)$. We deduce that
\[\pl_j\log\,\det(I-K)=-{\rm Tr}\;(I-K)\inv\,\pl_jK=R_{jj}(\x_j,\,\x_j).\]
Hence $\D \log\,\det(I-K)={\rm Tr}\;r$, and so it follows from (\ref{eq3}) that
\[{\cal {D}}^2\,\log\,\det(I-K)=-{\rm Tr}\;q\,\on\,\tl q\]
since the trace of $[\t,\,r]$ equals zero. This gives the representation
\[\det(I-K)=\exp\left\{-\int_0^\iy\eta\,{\rm Tr}\;
q(\x+\eta)\,\on\,\tl q(\x+\eta)\;d\eta\right\}.\]
Here the determinant is evaluated at $(\x_1,\,\ldots,\,\x_m)$ and
in the integral $\x+\eta$ is shorthand for $(\x_1+\eta,\,\ldots,\,\x_m+\eta).$

If $m=1$ the commutators drop out, $q=\tl q$, equations (\ref{eq1}) and (\ref{eq2})
are Painlev\'e~II and these are the previously known results.

\newpage
\begin{center}{\bf III. Proof}\end{center}

The proof will follow along the lines of the derivation in \cite{Airy} for the case $m=1$.
There the kernel was ``integrable'' in the sense that its commutator with $M$, the
operator of multiplication by $x$, was of finite rank. The same was then true of the
resolvent kernel, which was useful. But now our kernel is not integrable, so there will necessarily be some 
differences.

With $D=d/dx$ we compute that
\[[D,\,K]_{ij}=-\A(x)\,\A(y)\,\ch_j(y)+L_{ij}(x,\x_j)\,\dl(y-\x_j)
+(\t_i-\t_j)\,K_{ij}(x,\,y).\]
Equivalently,
\[ [D,\,K]=-A(x)\,\on\,A(y)\,\ch(y)+L\,\dl+[\t,\,K],\]
where $\dl=\sum_j\dl_j$, multiplication by the matrix ${\rm diag}\,(\dl(x-\x_j))$,
and $L$ is the operator with kernel $L_{ij}(x,y)$.
(For clarity we sometimes write the kernel of an operator
in place of the operator itself.)
To obtain $[D,\,R]$ we replace $K$ by $K-I$ in the commutators and left- and 
right-multiply by $\r=(I-K)\inv$. The result is
\be[D,\,R]=-Q(x)\,\on\,\tl Q(y)+R\,\dl\,\r+[\t,\,\r].
\footnote{Because of the fact $\r\,L\,\ch=R$ and our interpretation of $R_{ij}(x,\,\x_j)$ as 
$R_{ij}(x,\,\x_j+)$
we are able to write $R\,\dl\,\r$ in place of $\r\,L\,\dl\,\r$.}
\label{DRcom}\ee

We have already defined the matrix functions $Q$ and $\tl Q$  and we define 
\[P=(I-K)\inv A',\ \ \ u=(\tl Q,\,\A)=\int \tl Q(x)\,\A(x)\,dx.\]
It follows from (\ref{DRcom}) and the fact that
$\t$ and $A$ commute that
\be Q'=P-Q\,\on\,u+R\,\dl\,Q+[\t,\,Q].\footnote{The meaning of $\dl$ here and later is this: 
If $U$ 
and $V$ are matrix functions then $U\,\dl\,V$ is the matrix
with $i,j$ entry $\sum_kU_{ik}(\x_k)\,V_{kj}(\x_k)$. Thus $R\,\dl\,Q$ is the matrix
function with $i,j$ entry $\sum_kR_{ik}(x,\,\x_k)\,Q_{kj}(\x_k)$. This makes it compatible 
with our use of $\dl$ also as a multiplication operator so that, for example, 
$(R\,\dl\,\r)\,(A)=R\,\dl\,(\r\,A)$.}\label{Q'}\ee

Next, it follows from (\ref{plK}) that
\be\pl_j\,R=-R\,\dl_j \,\rho,\label{plR}\ee
and it follows from this that 
$\pl_j Q=-R\,\dl_j\,Q.$
Summing over $j$, adding to (\ref{Q'}) and evaluating at $\x_k$ give
\[\D Q(\x_k)=P(\x_k)-Q(\x_k)\,\on\,u+[\t,\,Q(\x_k)].\]
If we define $p_{ij}=P_{ij}(\x_i)$
then we obtain
\be\D q=p-q\,\on\,u+[\t,\,q].\label{Dq}\ee

Next we use the facts that $D^2-M$ commutes with $L$ and that $M$ commutes with $\ch$. 
It follows that 
\[ [D^2-M,\,K]=[D^2-M,\,L\,\ch]=L\,[D^2-M,\,\ch]=L\,[D^2,\,\ch]=L\,(\dl\, D+D\,\dl).\]
It follows from this that
\[ [D^2-M,\,\r]=\r\,L\,\dl\,D\,\r+\r\,L\,D\,\dl\,\r.\]
Applying both sides to $A$ and using the fact that $(D^2-M)A=0$ we obtain
\be Q''(x)-x\,Q(x)=\r\,L\,\dl\,Q'+\r\,L\,D\,\dl\,Q.\label{Q''}\ee
The first term on the right equals $R\,\dl\,Q'$. For the second term observe that 
\[\r\,L\,D\,\ch=\r\,L\,\ch\,D+\r\,L\,[D,\,\ch]=R\,D+\r\,L\,\dl,\]
so we can interpret that term as $-R_y\,\dl\,Q$ (the subscript denotes partial derivative) 
where $-R_y(x,y)$ is 
interpreted as not containing the delta-function summand which arises from the jumps of $R$.  With 
this interpretation of $R_y$
we can write the second term on the right as $-R_y\,\dl\,Q$. Thus,
\[Q''(x)-x\,Q(x)=R\,\dl\,Q'-R_y\,\dl\,Q.\]

Using this we obtain from (\ref{Q'})
\[P'=x\,Q(x)+R\,\dl\,Q'-R_y\,\dl\,Q+Q'\,\on\,u-R_x\,\dl\,Q-[\t,\,Q'],\]
and then from (\ref{Q'}) once more
\[P'=x\,Q(x)+R\,\dl\,(P-Q\,\on\,u+R\,\dl\,Q+[\t,\,Q])-R_y\,\dl\,Q\]
\[+(P-Q\,\on\,u+R\,\dl\,Q+[\t,\,Q])\,\on\,u-R_x\,\dl\,Q-[\t,\,P-Q\,\on\,u+R\,\dl\,Q+[\t,\,Q]\,].\]
It follows from (\ref{DRcom}) that
\[R_x+R_y=-Q(x)\,\on\,\tl Q(y)+R\,\dl\,R+[\t,\,\r].\]
(We replaced
$R\,\dl\,\r$ by $R\,\dl\,R$ since, recall, $R_y$ does not contain delta-function summands.)
We use this and also the identity $R\dl[\t,Q]-[\t,R\dl Q]=-[\t,R\dl] Q$, and the fact that
$\dl$ and $\t$ commute. The result is that \pagebreak
\[P'=x\,Q(x)+R\,\dl\,P+Q(x)\,\on\,\tl Q\dl Q+(P-Q\,\on\,u+[\t,\,Q])\,\on\,u\]
\[ -2[\t,\,R]\,\dl\,Q-[\t,\ P-Q\,\on\, u+[\t,\,Q]\,].\]

It follows from (\ref{plR}) that $\pl_j P=-R\,\dl_j\,P$. Summing over $j$, 
adding to the above and evaluating at $\x_k$ give 
\[\D P(\x_k)=\x_k\,Q(\x_k)+Q(\x_k)\,\on\,\tl Q\dl Q+(P(\x_k)-Q(\x_k)\,\on\,u+
[\t,\,Q(\x_k)])\,\on\,u\]
\[-2\,[\t,\,R(\x_k,\,\cdot\,)]\,\dl\,Q-[\t,\ P(\x_k)-Q(\x_k)\,\on\,u+[\t,\,Q(\x_k)]\,].\]
Hence $\D p$ is equal to
\[\x\,q+q\,\on\,\tl q\,q+(p-q\,\on\,u+[\t,\,q])\,\on\,u
-2\,[\t,\,r]\,q-[\t,\, p-q\,\on\,u+[\t,\,q]\,].\]
Equivalently, in view of (\ref{Dq}),
\be \D p=\x\,q+q\,\on\,\tl q\,q+\D q\,\on\, u
-2\,[\t,\,r]\,q-[\t,\, \D q].\label{Dp}\ee

Let us compute $\D u$. We have
\[u_{ij}=\int\int\A(x)\,\ch_i(x)\,\r_{ij}(x,y)\,\A(y)\,dx\,dy,\]
and so
\[\pl_k\,u_{ij}=-\dl_{ik}\,\int\A(\x_k)\,\r_{kj}(\x_k,y)\,\A(y)\,dy\]
\[-\int\int\A(x)\,\ch_i(x)\,[R_{ik}(x,\x_k)\,\r_{kj}(\x_k,y)]\,\A(y)\,dx\,dy,\]
where we use (\ref{plR}) again.
This is equal to
\[-\dl_{ik}\,\A(\x_k)\,Q_{kj}(\x_k)-\Big(\tl Q_{ik}(\x_k)-\dl_{ik}\,\A(\x_k)\Big)
\,Q_{kj}(\x_k),\]
and so
\be\pl_k\,u_{ij}=-\tl Q_{ik}(\x_k)\,Q_{kj}(\x_k).\label{plu}\ee
This gives
\be\D u=-\tl q\,q.\label{Du}\ee

Next, we find from (\ref{plR}) and (\ref{DRcom}) that
\[\D R(\x_j,\,\x_k)=-Q(\x_j)\,\on\,\tl Q(\x_k)
+[\t,\,R(\x_j,\,\x_k)].\]
This gives
$\D r=-q\,\on\,\tl q+[\t,\,r]$, which is equation (\ref{eq3}).

To get equation (\ref{eq1}) we apply $\D$ to (\ref{Dq})
and use (\ref{Dp}) and (\ref{Du}). We find that 
\[{\cal {D}}^2\,q=\x\,q+q\,\on\,\tl q\,q+\D q\,\on\, u
-2\,[\t,\,r]\,q
-[\t,\,\D q]-\D q\,\on\, u+q\,\on\,\tl q\,q+[\t,\,\D q]\]
\[=\x\,q+2\,q\,\on\,\tl q\,q-2\,[\t,\,r]\,q,\]
which is (\ref{eq1}).

Finally, to get equation (\ref{eq2}) we use the fact
that $\ch_j(y)\,\r_{jk}(y,x)$ is equal
to $\ch_k(x)$ times $\r_{kj}'(x,y)$, where $\r'$ is the resolvent kernel
for the matrix kernel with $i,j$ entry $L_{ji}(x,y)\,\ch_j(y)$. Hence
$\tl Q_{jk}(x)$ is equal to $\ch_k(x)$ times the $Q_{kj}(x)$ associated with $L_{ji}$.
Consequently for all the differentiation formulas we have for the $Q_{kj}(\x_k)$, etc., 
there are analogous formulas for the $\tl Q_{jk}(\x_k)$, etc.. The difference is
that we have to reverse subscripts and replace $r$ by $r^t$
and $\t$ by $-\t$. The
upshot is that, by computations analogous to those used to derive (\ref{eq1}), we derive
another equation which can be obtained from (\ref{eq1})
by making the replacements $q\to\tl q^t$, $\tl q\to q^t$,  $r \to r^t$,
$\t\to-\t$ and then taking transposes. The result is equation~(\ref{eq2}).\sp

\begin{center}{\bf Acknowledgment}\end{center}

This work was supported by the National Science Foundation through grants
DMS-9802122 and DMS-9732687.\sp


\begin{thebibliography}{9999}

\bibitem{adler} M.~Adler and P.~van Moerbeke, \textit{A PDE for the
joint distributions of the Airy process}, preprint,
arXiv:\,math.PR/0302329.

\bibitem{AD} D.~Aldous and P.~Diaconis, \textit{Longest increasing
subsequences: from patience sorting  to the Baik-Deift-Johansson
theorem}, Bull.\ Amer.\ Math.\ Soc.\ \textbf{36} (1999), 413--432.


\bibitem{BDJ} J.~Baik, P.~Deift and K.~Johansson,
\textit{On the distribution of the length of the longest increasing subsequence
in a random permutation}, J.\ Amer.\ Math.\ Soc.\ \textbf{12} (1999), 1119--1178.

\bibitem{J} K.~Johansson, \textit{Discrete polynuclear growth and determinantal processes},
Comm.\ Math.\ Phys.\ \textbf{242} (2003), 277--329.

\bibitem{KPZ} M.~Kardar, G.~Parisi and Y.~Z.~Zhang, 
\textit{Dynamic scaling of growing interfaces}, Phys.\ Rev.\ Letts.\
\textbf{56} (1986), 889--892.

\bibitem{PS} M.~Pr\"ahofer and H.~Spohn, \textit{Scale invariance of the PNG droplet and
the Airy process}, J.\  Stat.\  Phys.\  {\bf 108} (2002), 1071--1106.

\bibitem{Airy} C.~A.~Tracy and H.~Widom, \textit{Level-spacing distributions and the Airy
kernel}, Comm.\  Math.\  Phys.\  {\bf 159} (1994), 151--174.

\end{thebibliography}
\end{document}